\documentstyle[12pt]{article}

\input amssym.def
\input amssym

\newif\ifUS
\ifUS

\else

\fi

\def\twelveamsfonts{
 \font\twelvemsa=msam10 scaled 1200
 \font\eightmsa=msam8
 \font\sixmsa=msam6
 \font\twelvemsb=msbm10 scaled 1200
 \font\eightmsb=msbm8
 \font\sixmsb=msbm6
 \font\twelvembi=cmmib10 scaled 1200
 \font\eightmbi=cmmib8
 \font\sixmbi=cmmib6
 \textfont\msafam\twelvemsa
 \scriptfont\msafam\eightmsa
 \scriptscriptfont\msafam\sixmsa
 \textfont\msbfam\twelvemsb
 \scriptfont\msbfam\eightmsb
 \scriptscriptfont\msbfam\sixmsb}

\newif\iftwelve
\font\rmten=cmr10

\setbox0\hbox{\rm x} \setbox2\hbox{\rmten x}
\ifnum\ht0>\ht2 \twelvetrue\fi
\iftwelve\twelveamsfonts\fi

\let\bls\baselineskip \let\nt\noindent
\let\vp\vphantom 
\let\tsize\textstyle  
\def\vsk#1>{\vskip#1\bls} \def\vv#1>{\vadjust{\vsk#1>}}
\def\,{\relax\ifmmode\mskip\thinmuskip\relax\else\kern.16667em\fi}
\def\;{\relax\ifmmode\mskip\thickmuskip\relax\else\kern.27777em\fi}
\def\!{\relax\ifmmode\mskip-\thinmuskip\relax\else\kern-.16667em\fi}
\def\KZ/{{\sl KZ\/}} \def\qKZ/{{\sl qKZ\/}} \def\qKZB/{{\sl qKZB\/}}
\def\KZv/{Knizh\-nik-Za\-mo\-lod\-chi\-kov}
\def\&{.\kern.1em} \def\itl#1{{\it #1\/}} \def\qqq{\qquad\quad}
\def\ftext#1{{\let\thefootnote\relax\footnotetext{\vsk-.8>\noindent #1}}}

\let\Sum\sum \def\sum{\Sum\limits}
\let\Prod\prod \def\prod{\Prod\limits}
\def\Plus{\bigoplus\limits}
\def\diag{\mathop{\rm diag\;\!}}
\def\End{\mathop{\rm End\;\!}}
\def\Lim{\mathop{\rm lim\;\!}\nolimits}
\def\Re{\mathop{\rm Re\;\!}}
\def\Im{\mathop{\rm Im\;\!}}
\def\sgn{\mathop{\rm sgn\;\!}}

\let\eps\epsilon
\let\la\lambda
\let\om\omega

\def\C{{\Bbb C}}
\def\Z{{\Bbb Z}}

\let\ge\geqslant
\let\geq\geqslant
\let\le\leqslant
\let\leq\leqslant

\def\jmax{j_{\rm max}}
\def\mp{m^*}
\def\bfsis{\poorbold{$\scriptstyle\sigma$}}
\def\bfsi{\poorbold{$\sigma$}}

\def\sh{\mathop{\rm sh}}

\def\bea{\begin{eqnarray}}
\def\ena{\end{eqnarray}}
\def\beq{\begin{equation}}
\def\eeq{\end{equation}}
\def\no{\nonumber}
\def\nus{\nu_{1}, \cdots , \nu_{n-1}}
\def\lambdas{\lambda_{1}, \cdots , \lambda_{n}}
\def\qed{\quad$\square$}

\def\proof{\noindent{\it Proof.}\quad}
\newtheorem{thm}{Theorem}[section]
\newtheorem{prop}[thm]{Proposition}
\newtheorem{coro}[thm]{Corollary}

\makeatletter
\newbox\p@b@ld
\def\poorbold#1{\setbox\p@b@ld\hbox{#1}\kern-.01em\copy\p@b@ld\kern-\wd\p@b@ld
 \kern.02em\copy\p@b@ld\kern-\wd\p@b@ld\kern-.012em\raise.02em\box\p@b@ld}

\@addtoreset{thm}{section}
\@addtoreset{equation}{section}
\makeatother

\begin{document}

\begin{center}
\vp1
{\Large \bf Determinant Formula for Solutions of
\vsk.35>
the $U_q(sl_n)$ \qKZ/ Equation at $|q|=1$}
\vsk2>
{Tetsuji Miwa$^{\,\star}$, \;Yoshihiro Takeyama$^{\,\star\,\diamond}$,
\;Vitaly Tarasov$^{\,*}$}
\ftext{$^{\diamond\,}$Research Fellow of the Japan Society for
the Promotion of Science.}
\vsk1.5>
{\it $^\star$Research Institute for Mathematical Sciences,
Kyoto University, Kyoto 606, Japan
\vsk.5>
$^*$St\&Petersburg Branch of Steklov Mathematical Institute\\[2pt]
Fontanka 27, St\&Petersburg \,191011, Russia}
\end{center}
\vsk1.75>

{\narrower\nt
{\bf Abstract.}\enspace
We construct the hypergeometric solutions for the quantized \KZv/ equation with
values in a tensor product of vector representations of $U_q(sl_n)$ at $|q|=1$
and give an explicit formula for the corresponding determinant in terms of
the double sine function.
\vsk1.4>}
\vsk0>
\thispagestyle{empty}

\section*{Introduction}
In this paper we study the hypergeometric solutions of the quantized \KZv/
(\qKZ/) equation with values in a tensor product of vector representations
of $U_q(sl_n)$, see Section 1 for the precise formulation of the problem.
It is known that the \qKZ/ equation respects the weight decomposition of
the tensor product. For each weight subspace we construct a fundamental matrix
solution of the \qKZ/ equation and explicitly calculate the corresponding
determinant, see Theorem \ref{thm6.1}.

Formal integral representations for solutions of the \qKZ/ equation in
the $sl_n$ case, both in the rational and trigonometric situation, were
constructed in \cite{TV1}. Though to write down the phase function explicitly
in the trigonometric situation it had been assumed in \cite{TV1} that
the multiplicative step $p$ of the \qKZ/ equation is inside the unit circle:
$0<|p|<1$, all the construction in \cite{TV1} used only difference equations
for the phase function and after obvious modifications remained valid for
an arbitrary step $p\ne 0,1$. However, the problem of integrating the formal
integral representations suitably and getting in this way actual solutions of
the \qKZ/ equation is much more analytically involved; one can see this looking
at the $sl_2$ case.

In the last four years the hypergeometric solutions of the \qKZ/ equation
in the $sl_2$ case were studied quite well. The generic situation was
considered in \cite{TV2} (the rational case) and in \cite{TV3}
(the trigonometric case for ${0<|p|<1}\,$). The construction was generalized
to the trigonometric case for $|p|=1$ in \cite{MT1} and to the elliptic case
of the quantized \KZv/-Bernard (\qKZB/) equation in \cite{FTV1}, \cite{FTV2}.

If some of the representations are finite-dimensional, the situation is no
more generic. Rather detailed study of this case has been done in \cite{MV1};
see also \cite{S}, \cite{JKMQ}, \cite{NPT}, \cite{T1} for some important
particular cases.

Less is known for the case of ${n>2}$. Some integral formulae for solutions of
the \qKZ/ equation were obtained in \cite{S}, \cite{KQ}, \cite{N}, \cite{MT2},
but in all the considered cases solutions takes values in a tensor product of
vector representations. Recently Varchenko and the third author have managed to
extend the costruction of \cite{TV2}, \cite{TV3} to the higher rank case and
get solutions taking values in a tensor product of arbitrary highest weight
representations \cite{TV4}. Let us also mention a paper \cite{M}, where
integral formulae for solutions of another type of the \qKZ/ equation were
suggested.

In this paper we evaluate a determinant of a certain matrix whose entries are
given by multidimensional integrals of $q$-hypergeometric type. In the case
of ordinary multidimensional hypergeometric integrals a problem of evaluating
similar determinants appears, say, in studying arrangements of hyperplanes,
and several results have been obtained in this direction, see for instance
\cite{V1}, \cite{L}, \cite{LS}, \cite{DT}, \cite{MTV}, \cite{MV2}. In some
particular cases these determinant formulae have another meaning; namely they
imply that under certain assumptions the hypergeometric solutions of the
differential \KZv/ equation form a basis of solutions \cite{SV}, \cite
{V2}. There are similar determinant formulae for solutions of the \qKZ/
equation in the $sl_2$ case.  They have been obtained for the rational case in
\cite{TV2}, \cite {T1}, and for the trigonometric case in \cite{TV3} for
${0<|p|<1}$ and in \cite{MT1} for $|p|=1$. It turns out that there is a nice
connection of constructions given in \cite{TV3} and \cite {MT1}, which in
particularly allows to derive the determinant formula for $|p|=1$ from the
determinant formula for ${0<|p|<1}$. This subject will be addressed elsewhere
\cite{T2}.

The paper is organized as follows. The first section contains preliminaries
and precise definitions on the \qKZ/ equation. In Section 2 we construct the
hypergeometric pairing and give integral formulae for solutions of the \qKZ/
equation. The main result of the paper is formulated in Section 3, see Theorem
\ref{thm6.1}. We show that both the left hand side and the right hand side of
formula (\ref{107}) satisfy the same system of difference equations and have to
be proportional. To compute the proportionality coefficient we study suitable
asymptotics of the hypergeometric solutions. We see that the proportionality
coefficient splits into a product of contributions of each tensor factor, which
are calculated in Section 5. In the last Section we complete the proof of
Theorem \ref{thm6.1}. A short Appendix contains the necessary information of
the double sine function for the convenience of the reader.

\section*{Acknoledgement}
V\&Tarasov thanks the RIMS, Kyoto University for warm hospitality during his
stay there, when most of the results of this paper had been obtained.
He also thanks A\&Varchenko for valuable remarks.

\section{The quantized \KZv/ equation}
Consider the vector representation $V$ of $sl_n$:
$$
\tsize V\,=\,\Plus_{j=0}^{n-1}\,\C\,v_j\,.
$$
Let $\eps_1,\ldots,\eps_n$ be the fundamental weights of $sl_n$. We consider
the basis vectors $v_0,\ldots v_{n-1}$ as weight vectors with respect to
the Cartan subalgebra of $sl_n$ with weights $\eps_1,\ldots,\eps_n$,
respectively.

Fix a complex number $\rho$ and a weight $\bar\mu=\sum_{j=1}^n\mu_j\eps_j$.
Consider a diagonal matrix
$$
D(\bar\mu)\,=\,\diag(e^{2\pi i\mu_1},\ldots,e^{2\pi i\mu_n})
$$
and a matrix $R^{(\rho)}(\beta)\in\End(V\otimes V)$ with the following entries:
$R^{(\rho)}(\beta)_{jj}^{jj}=1$,
\bea
R^{(\rho)}(\beta)_{jk}^{jk}\,=\,
{\sh{\pi\over\rho}\beta\over\sh{\pi\over\rho}(\beta-{2\pi i\over n})} \label{eq01}
\ena
for $j\ne k$,
\beq
R^{(\rho)}(\beta)_{kj}^{jk}\,=\,
-\;{e^{{\pi\over\rho}\beta}\,\sh{2\pi^2i\over\rho n}\over
\sh{\pi\over\rho}(\beta-{2\pi i\over n})}\;,\qqq
R^{(\rho)}(\beta)_{jk}^{kj}\,=\,
-\;{e^{-{\pi\over\rho}\beta}\,\sh{2\pi^2i\over\rho n}\over
\sh{\pi\over\rho}(\beta-{2\pi i\over n})}
\label{eq02}
\eeq
for $j<k$, and $R^{(\rho)}(\beta)_{lm}^{jk}=0$, otherwise. We have
$$
R^{(\rho)}(\beta)\,v_l\otimes v_m\,=\,\sum_{j,k=0}^{n-1}
R^{(\rho)}(\beta)^{lm}_{jk}\,v_j\otimes v_k\,.
$$

Fix a complex number $\la$ and define the \qKZ/ operators
$K_1^{(\rho)},\ldots K_N^{(\rho)}$ acting in the tensor product
$V^{\otimes N}$:
\bea
K^{(\rho)}_m(\beta_1,\ldots,\beta_N)\,=\,
R^{(\rho)}_{m,m-1}(\beta_m-\beta_{m-1}-\la i)
\ldots R^{(\rho)}_{m1}(\beta_m-\beta_1-\la i)
\no
\\[4pt]
{}\times\,D_m(-\bar\mu/\rho)\,R^{(\rho)}_{mN}(\beta_m-\beta_N)\ldots
R^{(\rho)}_{m,m+1}(\beta_m-\beta_{m+1})\,.\kern-17pt
\ena
In this paper we consider the \qKZ/ equation for a function
$f(\beta_1,\ldots,\beta_N)$ taking values in $V^{\otimes N}$, which is
the following system of difference equations:
\beq
f(\beta_1,\ldots \beta_m-\la i,\ldots \beta_N)\,=\,
K^{(\rho)}_m(\beta_1,\ldots,\beta_N)\,f(\beta_1,\ldots,\beta_N)\,,
\quad m=1,\ldots N.\kern-2em
\label{qKZ}
\eeq
We also consider the mirror \qKZ/ equation for a similar function
$g(\beta_1,\ldots,\beta_N)$:
\beq
g(\beta_1,\ldots \beta_m-\rho i,\ldots \beta_N)\,=\,
K^{(\la)}_m(\beta_1,\ldots,\beta_N)\,g(\beta_1,\ldots,\beta_N)\,,
\quad m=1,\ldots N.\kern-2em
\label{qKZm}
\eeq
The \qKZ/ operators respects the $sl_n$ weight decomposition of the tensor
product. Therefore, one can consider solutions of the \qKZ/ and mirror \qKZ/
equations taking values in a weight subspace $(V^{\otimes N})_\xi$ for any
given weight $\xi$.

In this paper for any given weight $\xi$ such that $(V^{\otimes N})_\xi$
is nontrivial we will construct a function
$\Psi_\xi(\beta_1,\ldots,\beta_N)$ taking values in
$(V^{\otimes N})_\xi\otimes(V^{\otimes N})_\xi$, which solves the \qKZ/
equation (\ref{qKZ}) in the first tensor factor and solves the mirror \qKZ/
equation (\ref{qKZm}) in the second tensor factor. Also we will compute the
determinant $\det\Psi_\xi(\beta_1,\ldots,\beta_N)$.

The matrix $R^{(\rho)}(\beta)$ is the $R$-matrix associated with the tensor
product of the evaluation vector representations of $U_q(\widehat{sl_n})$ for
$$
q\,=\,e^{-{2\pi^2 i\over \rho n}}\,.
$$
Similarly, $R^{(\la)}(\beta)$ is associated with $U_{q'}(\widehat{sl_n})$ for
$$
q'\,=\,e^{-{2\pi^2 i\over\la n}}\,.
$$
All over this paper we assume that $\rho$ and $\la$ are real positive.
Thus, under our assumptions we have that
$$
|q|\,=\,|q'|\,=\,1\,.
$$
However, it is clear from the consideration that all our construction remain
valid if $\rho$ and $\la$ have small enough imaginary parts of arbitrary sign.
Therefore, $q$ and $q'$ can deviate from the unit circle and vary in a narrow
annulus.

In addition to the reality and positivity of $\rho$ and $\la$ we assume that
both of them and their ratio are not rational. Sometimes we take $\rho$ and
$\la$ to be sufficiently large.
\section{Integral formulae for solutions}

For non-negative integers $\nu_{1}, \cdots , \nu_{n-1}$ satisfying
\bea
N=\nu_{0} \ge \nu_{1} \ge \cdots \ge \nu_{n-1} \ge \nu_{n}=0,
\ena
we denote by ${\cal Z}_{\nus}$ the set of all $N$-tuples
$J=(J_{1}, \cdots , J_{N}) \in \left( \Z_{\ge 0} \right)^{N}$
such that
\bea \# \{ r\,;\,J_{r} \ge j \}=\nu_{j}.
\ena
For $J=(J_{1}, \cdots , J_{N}) \in {\cal Z}_{\nus}$, we set
\bea
{\cal N}_{j}^{J}=\{ r\,;\,J_{r} \ge j \}
\ena
and define integers
$r_{j, m}^{J}\,, \ (\,0 \le j \le n-1\,,\ \ 1 \le m \le \nu_{j}\,)$ as follows:
\bea
{\cal N}_{j}^{J}=\{ r_{j, 1}^{J}, \cdots , r_{j, \nu_{j}}^{J} \},
\quad r_{j, 1}^{J}< \cdots < r_{j, \nu_{j}}^{J}.
\ena
We have, in particular, $r_{0, m}^{J}=m$.

Now we set
\bea
& & w_{J}^{(\rho)}(\{ \gamma_{j, m} \} ; \beta_{1}, \cdots , \beta_{N})=
{\rm Skew}_{n-1} \circ \cdots \circ {\rm Skew}_{1}
g_{J}^{(\rho)}(\{\gamma_{j, m}\}; \beta_{1}, \cdots , \beta_{N}),
\\[6pt]
& & g_{J}^{(\rho)}(\{\gamma_{j, m}\}; \beta_{1}, \cdots , \beta_{N})\no
\\[2pt]
& & =\prod_{j=1}^{n-1} \Bigl\{
\prod_{1 \le m<m' \le \nu_{j}}
\sh\frac{\pi}{\rho}(\gamma_{j, m'}-\gamma_{j, m}-\frac{2\pi i}{n})
\prod_{m=1}^{\nu_{j}} \Bigl\{
e^{-\frac{\pi}{\rho}(\gamma_{j, m}-\gamma_{j-1, \mp(J, j, m)})} \no \\
& & {}\times
\prod_{m' \atop r_{j-1, m'}^{J}<r_{j, m}^{J}}
\sh\frac{\pi}{\rho}(\gamma_{j, m}-\gamma_{j-1, m'}+\frac{\pi i}{n})\!\!\!\!\!
\prod_{m' \atop r_{j, m}^{J}<r_{j-1, m'}^{J}}
\sh\frac{\pi}{\rho}(\gamma_{j, m}-\gamma_{j-1, m'}-\frac{\pi i}{n})
\Bigr\} \Bigr\}.
\ena
The notation in the above formulae is as follows.
The operator ${\rm Skew}_{j}$ is the skew-symmetrization
with respect to the variables $\{ \gamma_{j, m}\}_{m=1, \cdots , \nu_{j}}$:
\bea
{\rm Skew}_{j}X(\gamma_{j, 1}, \cdots , \gamma_{j, \nu_{j}})=
\sum_{\sigma \in S_{\nu_{j}}}(\sgn\sigma)
X(\gamma_{j, \sigma(1)}, \cdots , \gamma_{j, \sigma(\nu_{j})}).
\ena
The integer $\mp(J,j,m)$ is uniquely determined by the condition
\bea
r_{j,m}^{J}=r_{j-1,\mp(J, j, m)}^{J},
\ena
and $\gamma_{0, m}=\beta_{m}$. We abbreviate
$w_{J}^{(\rho)}(\{\gamma_{j, m}\}; \beta_{1}, \cdots , \beta_{N})$ to
$w_{J}^{(\rho)}(\beta_{1}, \cdots , \beta_{N})$ when the dependence on
the abbreviated variables is irrelavant.

Set
\bea
{\cal F}_{\nus}^{(\rho)}(\beta_{1}, \cdots , \beta_{N})=
\sum_{J \in {\cal Z}_{\nus}}{\Bbb C} \,\,w_{J}^{(\rho)}(\beta_{1}, \cdots,
\beta_{N}).
\ena
In the following, we define a pairing between
\bea
w_{J}^{(\rho)} \in {\cal F}_{\nus}^{(\rho)}(\beta_{1}, \cdots ,
\beta_{N}) \quad
{\rm and} \quad
w_{J}^{(\lambda)} \in {\cal F}_{\nus}^{(\lambda)}(\beta_{1}, \cdots ,
\beta_{N}).
\ena
We use
\bea
\varphi(x)=\frac{1}{S_{2}(ix-\frac{\pi}{n})S_{2}(-ix-\frac{\pi}{n})}, \quad
\psi(x)=\frac{1}{S_{2}(ix+\frac{2\pi}{n})S_{2}(-ix+\frac{2\pi}{n})},
\ena
where $S_{2}(x)=S_{2}(x | \rho , \lambda )$ is the double sine function with
periods $\rho$ and $\lambda$.

For $J, J' \in {\cal Z}_{\nus}$, we set
\bea
& &
I(w_{J}^{(\rho)}, w_{J'}^{(\lambda)})=\left( \prod_{j=1}^{n-1}\prod_{m=1}^{\nu_{j}}
\int_{C_{j}} d\gamma_{j, m} \right)
K(\{\gamma_{j, m}\}; \beta_{1}, \cdots , \beta_{N}; \mu_{1}, \cdots , \mu_{n}) \no \\[4pt]
& & \quad{}\times
w_{J}^{(\rho)}(\{\gamma_{j, m}\}; \beta_{1}, \cdots , \beta_{N})
w_{J'}^{(\lambda)}(\{ \gamma_{j, m}\} ; \beta_{1}, \cdots , \beta_{N}). \label{001}
\ena
Here the function $K$ is defined by
\bea
&&K(\{\gamma_{j,m}\};\beta_{1},\cdots,\beta_{N};\mu_{1},\cdots,\mu_{n})=
e^{\frac{2\pi}{\rho \lambda}\sum_{j=0}^{n-1}\sum_{m=1}^{\nu_{j}}
\gamma_{j,m}(\mu_{j+1}-\mu_{j})}\no\\[6pt]
&&{}\times\prod_{j=1}^{n-1}\left\{\prod_{m=1}^{\nu_{j}}
\prod_{m'=1}^{\nu_{j-1}}\varphi(\gamma_{j,m}-\gamma_{j-1,m'})
\prod_{1\le m<m'\le\nu_{j}}\psi(\gamma_{j,m}-\gamma_{j,m'})\right\},
\label{Kdef}
\ena
where $\mu_{0}=0$. We say the variable $\gamma_{j,m}$ belongs to the point
$r_{j,m}$. Then, $J_r$ is the number of integral variables which belong
to the point $r$.

The contour $C_{j}$ for $\gamma_{j,m}, \,(m=1, \cdots , \nu_{j})$
is a deformation of the real line $(-\infty, \infty)$
such that the poles at
\bea
\gamma_{j-1, m'}-\frac{\pi i}{n}+\rho i \Z_{\ge 0}+\lambda i \Z_{\ge 0}, \quad
\gamma_{j, m'}+\frac{2\pi i}{n}+\rho i \Z_{\ge 0}+\lambda i \Z_{\ge 0}
\label{CONT1}
\ena
are above $C_{j}$ and the poles at
\bea
\gamma_{j-1, m'}+\frac{\pi i}{n}-\rho i \Z_{\ge 0}-\lambda i \Z_{\ge 0}, \quad
\gamma_{j, m'}-\frac{2\pi i}{n}-\rho i \Z_{\ge 0}-\lambda i \Z_{\ge 0}
\label{CONT2}
\ena
are below $C_{j}$, where $\gamma_{0, m}=\beta_{m}$.

These conditions are not compatible if all the poles really exist.
Pinching of the integration contours by poles occurs for each triple
of variables $\gamma_{j,m_1},\gamma_{j,m_2},\gamma_{j-1,m}$. However,
we can improve the definition (\ref{001}) as follows. We have that
$$
I(w_{J}^{(\rho)}, w_{J'}^{(\lambda)})=
\sum_{{\bfsis}, {\bfsis'}\in\,
{\bf S}_{\nu_{1}} \times \cdots \times {\bf S}_{\nu_{n-1}}}\!\!
(\sgn{\bfsi})\,(\sgn{\bfsi'})\, \left( \prod_{j=1}^{n-1}\prod_{m=1}^{\nu_{j}}
\int_{C_{j}} d\gamma_{j, m} \right)
F_{J,J'\!,\,\bfsis,\,\bfsis'}(\{\gamma_{j,m}\})
$$
where ${\bfsi}=(\sigma_{1},\cdots,\sigma_{n-1})\in{\bf S}_{\nu_{1}}
\times\cdots\times{\bf S}_{\nu_{n-1}}\,$,
\ $\sgn{\bfsi}=\prod_{j=1}^{n-1}(\sgn\sigma_{j})\,$ and
\bea
F_{J,J'\!,\,\bfsis,\,\bfsis'}(\{\gamma_{j,m}\}) & \!\!\!=\!\!\! &
K(\{\gamma_{j, m}\};\beta_{1},\cdots,\beta_{N}; \mu_{1},\cdots,\mu_{n})
\no
\\[4pt]
& \!\!\!\times\!\!\! &
g_{J}^{(\rho)}(\{\gamma_{j,\sigma_{j}(m)}\};\beta_{1},\cdots,\beta_{N})\,
g_{J'}^{(\lambda)}(\{\gamma_{j,\sigma_{j}'(m)}\};\beta_{1},\cdots,\beta_{N})\,.
\no
\\[-.9\bls]\no
\ena
A partial integrand $F_{J,J'\!,\,\bfsis,\,\bfsis'}(\{\gamma_{j,m}\})$
does not have poles at some points (\ref{CONT1}) and (\ref{CONT2})
because of the zeros of
$g_{J}^{(\rho)}(\{\gamma_{j, \sigma(m)}\};\beta_{1},\cdots,\beta_{N})$.
So, given $J$ and $\bfsi$ there is a choice of integration contours
$C_j^{(J,\bfsis)}$ satisfying the required conditions for the actual poles of
$F_{J,J'\!,\,\bfsis,\,\bfsis'}(\{\gamma_{j,m}\})$ for arbitrary
$J'\!,\,\bfsi'$.
Similarly, given $J'$ and $\bfsi'$ there is a choice of integration contours
$C_{j,(J'\!,\bfsis')}$ satisfying the required conditions for the actual poles
of $F_{J,J'\!,\,\bfsis,\,\bfsis'}(\{\gamma_{j,m}\})$ for arbitrary $J,\,\bfsi$.
Finally, one can easily check that the integrals of the term
$F_{J,J'\!,\,\bfsis,\,\bfsis'}(\{\gamma_{j,m}\})$ over the contours
$C_j^{(J,\bfsis)}$ and over the contours $C_{j,(J'\!,\bfsis')}$ are equal.

\vsk.5>
In this paper we assume that $\rho$ and $\lambda$ are large positive.
Then, as we will see in the next section, there is a region of the parameters
$\mu_{1}, \cdots , \mu_{n}$ where the integral (\ref{001}) is absolutely
convergent (see (\ref{B})).
\vsk.5>

Consider the vector representation $V$ of $\,sl_n\,$:
\beq
\tsize V\,=\,\Plus_{j=0}^{n-1}\,\C\,v_j\,.
\eeq

\begin{thm}\label{thm2.2}
For $J \in {\cal Z}_{\nus}$, we set
\bea
f_{J}=\sum_{J'\in{\cal Z}_{\nus}}I(w_{J'}^{(\rho)}, w_{J}^{(\lambda)})v_{J'},
\quad {\it where} \quad v_{J'}=v_{J_{1}'}\otimes\cdots\otimes v_{J_{N}'}.
\ena
Then $f_{J}$ is a solution to (\ref{qKZ}).
\end{thm}

\proof
For $J=(J_{1}, \cdots , J_{N}) \in {\cal Z}_{\nus}$,
we set
$w_{J_{1}, \cdots , J_{N}}^{(\rho)}=w_{J}^{(\rho)}$.
In the same way as the proof of Lemma1 and Lemma 3 in \cite{MT2},
we can show the following formulae:
\bea
& & w_{J_{1}, \cdots , J_{k+1}, J_{k}, \cdots , J_{N}}^{(\rho)}
(\beta_{1}, \cdots , \beta_{k+1}, \beta_{k}, \cdots , \beta_{N})\no \\
& & {}=\sum_{J_{k}', J_{k+1}'}R^{(\rho)}(\beta_{k}-\beta_{k+1})_{J_{k}, J_{k+1}}^{J_{k}', J_{k+1}'}
w_{J_{1}, \cdots , J_{k}', J_{k+1}', \cdots , J_{N}}^{(\rho)}
(\beta_{1}, \cdots , \beta_{k}, \beta_{k+1}, \cdots , \beta_{N}), \label{pf00} \\
& & I(g_{J_{N}, J_{1}, \cdots , J_{N-1}}^{(\rho)}(\beta_{N}, \beta_{1}, \cdots , \beta_{N-1}),
w^{(\lambda)}(\beta_{1}, \cdots , \beta_{N}))|_{\beta_{N} \to \beta_{N}-\lambda i} \no \\
& & {}=e^{-\frac{2\pi i}{\rho}\mu_{J_{N}+1}}
I(g_{J_{1}, \cdots , J_{N}}^{(\rho)}(\beta_{1}, \cdots , \beta_{N-1}, \beta_{N}),
w^{(\lambda)}(\beta_{1}, \cdots , \beta_{N})), \label{pf01}
\ena
where $w^{(\lambda)} \in {\cal F}_{\nus}^{(\lambda)}(\beta_{1}, \cdots , \beta_{N})$ and
the left hand side of (\ref{pf01}) is understood
as the analytic continuation of the integral.
It is easy to prove Theorem \ref{thm2.2} from (\ref{pf00}) and (\ref{pf01}).
\qed
\vsk>

We note that the weight of the solution $\psi_{J}$ is given by
\bea
\sum_{j=1}^{n}\lambda_{j}\eps_{j},\quad{\rm where}\quad\lambda_{j}=\nu_{j-1}-\nu_{j}.
\ena
Now we set
\bea
\Psi_\xi(\beta_{1}, \cdots , \beta_{N})=
\sum_{J, J' \in {\cal Z}_{\nus}}I(w_{J}^{(\rho)}, w_{J'}^{(\lambda)})v_{J}^{(\rho)} \otimes v_{J'}^{(\lambda)},
\ena
where
\bea
\xi=\sum_{j=1}^{n}\lambda_{j}\eps_{j}.
\ena
Then $\Psi_\xi$ is the fundamental matrix solution mentioned in Section 1.

\section{Determinant formula for the solutions}

In the following sections we calculate the determinant
\bea
D_{\lambda_{1}, \cdots , \lambda_{n}}(\beta_{1}, \cdots , \beta_{N})=
\det{\left(I(w_{J}^{(\rho)}, w_{J'}^{(\lambda)})
\right)_{J, J' \in {\cal Z}_{\nus}}}.
\ena
The result is as follows.
\begin{thm}\label{thm6.1}
\bea
& &
\det{\left(I(w_{J}^{(\rho)}, w_{J'}^{(\lambda)})
\right)_{J, J' \in {\cal Z}_{\nus}}} \no \\
& & {}=
2^{-d_{\nus}\Lambda_{\lambdas}^{(0)}}
\exp{\Bigl( (N^{2}-N-2) \left(\frac{1}{\rho}+\frac{1}{\lambda}
\right) \frac{\pi^{2}i}{n}\;\Lambda_{\lambdas}^{(2)} \Bigr)}
\no \\
& & {}\times
\left(\frac{\sqrt{\rho\lambda}}{S_{2}(-\frac{2\pi}{n})} \right)^{\Lambda_{\lambdas}^{(0)}\sum_{j=1}^{n}(j-1)\lambda_{j}} \no \\
& & {}\times
\prod_{1 \le r' < r \le n}
\left\{ \prod_{a=0}^{\min\{\lambda_{r}-1, \lambda_{r'}\}}
\left( \frac{S_{2}(\mu_{r}-\mu_{r'}-\frac{\pi}{n}(\lambda_{r}+\lambda_{r'}-2a))}
{S_{2}(\mu_{r}-\mu_{r'}+\frac{\pi}{n}(\lambda_{r}+\lambda_{r'}-2a))} \right)
^{\lambda_{r}+\lambda_{r'} \atopwithdelims() a} \right\}
^{\Lambda_{\lambdas}^{(0)}\big/
{\lambda_{r}+\lambda_{r'} \atopwithdelims() \lambda_{r}}} \no \\
& & {}\times
\exp{\Bigl( \frac{2\pi}{\rho \lambda}
\left( \sum_{j=1}^{n}
{N-1 \atopwithdelims() \lambda_{1}, \cdots , \lambda_{j}-1, \cdots , \lambda_{n}}\mu_{j} \right)
\sum_{m=1}^{N}\beta_{m} \Bigr)} \no \\
& & \times{}
\left( \prod_{1 \le r < s \le N}\frac{S_{2}(i(\beta_{r}-\beta_{s})+
\frac{2\pi}{n})}{S_{2}(i(\beta_{r}-\beta_{s})-\frac{2\pi}{n})}
\right)^{\Lambda_{\lambdas}^{(2)}},
\label{107}
\ena
where
\bea
& & \lambda_{j}=\nu_{j-1}-\nu_{j}, \quad d_{\nus}=\sum_{j=1}^{n-1}(2\nu_{j}\nu_{j-1}+\nu_{j}^{2}-3\nu_{j}), \no \\
& & \Lambda_{\lambdas}^{(0)}={N \atopwithdelims() \lambdas}, \quad
\Lambda_{\lambdas}^{(2)}=
\sum_{1 \le j<k \le n}
{N-2 \atopwithdelims() \lambda_{1}, \cdots , \lambda_{j}-1, \cdots , \lambda_{k}-1, \cdots ,\lambda_{n}}. \no
\ena
\end{thm}

Note that
\bea
\Lambda_{\lambdas}^{(0)}\bigg/
{\lambda_{r}+\lambda_{r'} \atopwithdelims() \lambda_{r}}=
{N \atopwithdelims() \lambda_{1}, \cdots , \lambda_{r'-1}, \lambda_{r'+1},
\cdots , \lambda_{r-1}, \lambda_{r+1}, \cdots , \lambda_{n}, \lambda_{r'}+\lambda_{r}}
\ena
is a positive integer.

First, we determine the dependence on $\beta_{1}, \cdots , \beta_{N}$
of $D_{\lambdas}$.
{}From Theorem \ref{thm2.2}, we find that
\bea
\frac{D_{\lambdas}(\beta_{1}, \cdots , \beta_{m}-\lambda i, \cdots , \beta_{N})}
{D_{\lambdas}(\beta_{1}, \cdots , \beta_{m}, \cdots , \beta_{N})} &=&
\det_{\lambdas}K_{m}^{(\rho)}(\beta_{1}, \cdots , \beta_{N}), \label{201} \\
\frac{D_{\lambdas}(\beta_{1}, \cdots , \beta_{m}-\rho i, \cdots , \beta_{N})}
{D_{\lambdas}(\beta_{1}, \cdots , \beta_{m}, \cdots , \beta_{N})} &=&
\det_{\lambdas}K_{m}^{(\lambda)}(\beta_{1}, \cdots , \beta_{N}). \label{202}
\ena
Here $\det_{\lambdas}K_{m}$ stands for the determinant of the operator which
is a restriction of the operator $K_{m}$ to the weight subspace of the weight
$\sum_{j=1}^{n}\lambda_{j}\eps_{j}$.

Using formulae (\ref{eq01}) and (\ref{eq02}), we have
\bea
& & \det_{\lambdas}K_{m}^{(\rho)}(\beta_{1}, \cdots , \beta_{N})=\exp{\Bigl( -\frac{2\pi i}{\rho}\sum_{j=1}^{n}
{N-1 \atopwithdelims() \lambda_{1}, \cdots , \lambda_{j}-1, \cdots , \lambda_{n}}\mu_{j} \Bigr)} \no \\
& & {}\times
\left( \prod_{m'=1}^{m-1}\frac{\sh\frac{\pi}{\rho}(\beta_{m}-\beta_{m'}-\lambda i+\frac{2\pi i}{n})}
{\sh\frac{\pi}{\rho}(\beta_{m}-\beta_{m'}-\lambda i-\frac{2\pi i}{n})}
\prod_{m'=m+1}^{N}\frac{\sh\frac{\pi}{\rho}(\beta_{m}-\beta_{m'}+\frac{2\pi i}{n})}
{\sh\frac{\pi}{\rho}(\beta_{m}-\beta_{m'}-\frac{2\pi i}{n})} \right)
^{\Lambda_{\lambdas}^{(2)}}.
\ena

Now we set
\bea
E_{\lambdas}(\beta_{1}, \cdots , \beta_{N})&=&
\exp{\Bigl( \frac{2\pi}{\rho \lambda}
\left( \sum_{j=1}^{n}
{N-1 \atopwithdelims() \lambda_{1}, \cdots , \lambda_{j}-1, \cdots , \lambda_{n}}\mu_{j} \right)
\sum_{m=1}^{N}\beta_{m} \Bigr)} \no \\
& & \times{}
\left( \prod_{1 \le r < s \le N}\frac{S_{2}(i(\beta_{r}-\beta_{s})+\frac{2\pi}{n})}
{S_{2}(i(\beta_{r}-\beta_{s})-\frac{2\pi}{n})} \right)
^{\Lambda_{\lambdas}^{(2)}}.
\label{202.1}
\ena
Then by using (\ref{app3}) we can check that
$E_{\lambdas}(\beta_{1}, \cdots , \beta_{N})$ satisfies (\ref{201})
and (\ref{202})\,.  Therefore, we have

\begin{prop}\label{prop4.1}
\bea
D_{\lambdas}(\beta_{1}, \cdots , \beta_{N})=c_{\lambdas}(\mu_{1}, \cdots , \mu_{n} ; \rho, \lambda)
E_{\lambdas}(\beta_{1}, \cdots , \beta_{N}),
\label{202.2}
\ena
where $c_{\lambdas}(\mu_{1}, \cdots , \mu_{n} ; \rho, \lambda)$ is a constant
independent of $\beta_{1}, \cdots , \beta_{N}$.
\end{prop}

In order to determine $c_{\lambdas}(\mu_{1},\cdots ,\mu_{n};\rho,\lambda)$,
we consider the asymptotics of $D/E$ as
\bea
\beta_{1}\ll \cdots\ll \beta_{N}\,. \label{202.9}
\ena
This is in the next section.

\section{Asymptotics of the solutions}

First, we consider the asymptotics of $D_{\lambdas}$.

We denote the set of variables
\bea
\gamma_{j,m}\quad(0\leq j\leq n-1\,;\ \ 1\leq m\leq\nu_j)
\ena
by $\gamma$.
Fix a set of permutations $\sigma=(\sigma_1,\ldots,\sigma_{n-1})$;
$\sigma_j\in S_{\nu_j}$ $(1\leq j\leq n-1)$. We use $\sigma_0={\rm id}$.
We denote
\bea
\gamma_{j,\sigma_j(m)}\quad(0\leq j\leq n-1\,;\ \ 1\leq m\leq\nu_j)
\ena
by $\gamma_\sigma$.

Consider
\bea
F_{J,J',\sigma}(\gamma)
=K(\gamma)g_J(\gamma)g_{J'}(\gamma_\sigma),
\ena
where $K(\gamma)$ is given by (\ref{Kdef}).

In the following, we use the abbreviation $\beta_{ij}=\beta_{i}-\beta_{j}$.

\begin{prop}\label{prop1}
Suppose that $\beta_1<\cdots<\beta_n$ and $\gamma_{j,m}$'s are all real.
If $\lambda$ is sufficiently large, then
there exist positive constants $\varepsilon,C,\kappa$ independent
of the variables $\beta$ and $\gamma$ such that the following estimate holds.
\bea
|F_{J,J',\sigma}|&<&C
\exp{\Bigl( -\kappa\sum_{1\leq j\leq n-1\atop1\leq m\leq\nu_j}
|\gamma_{j,m}-\gamma_{j-1,\mp(J,j,m)}|\Bigr)}\nonumber\\
&\times&\exp{\Bigl( -{2\pi^2\over\rho\lambda n}
\sum_{1\leq r<s\leq N}(1-\delta_{J_r,J_s})\beta_{sr}
+{2\pi\over\rho\lambda}\sum_{r=1}^{N}\mu_{J_r+1}\beta_r \Bigr)}
\ena
if
\bea
{2\pi\over\rho\lambda}(\mu_{j+1}-\mu_j)>\varepsilon,\quad
{2\pi\over\rho\lambda}(\mu_n-\mu_1)<n\varepsilon.
\label{B}
\ena
\end{prop}

\proof
Throughout the proof, we set $r_{j, m}=r_{j, m}^{J}$.
We define new variables $\tilde\gamma$ by
\bea
\gamma_{j,m}=\tilde\gamma_{j,m}+\beta_{r_{j,m}}.
\ena
Note that $\tilde\gamma_{0,m}=0$. From (\ref{app2}), we have
\bea
e^{{2\pi\over\rho\lambda}\sum_{0\leq j\leq n-1\atop1\leq m\leq\nu_j}
(\mu_{j+1}-\mu_j)\gamma_{j,m}}
&=&
e^{{2\pi\over\rho\lambda}\sum_{1\leq j\leq n-1\atop1\leq m\leq\nu_j}
(\mu_{j+1}-\mu_j)\tilde\gamma_{j,m}
+{2\pi\over\rho\lambda}\sum_{1\leq r\leq N}\mu_{J_r+1}\beta_r},\nonumber\\
|\varphi(\gamma_{j,m}-\gamma_{j-1,m'})|&\leq&{\rm const.}\,
e^{-\pi({1\over\rho}+{1\over\lambda}+{2\pi\over\rho\lambda n})
|\tilde\gamma_{j,m}-\tilde\gamma_{j-1,m'}+\beta_{r_{j,m}r_{j-1,m'}}|},
\nonumber\\
|\psi(\gamma_{j,m}-\gamma_{j,m'})|&\leq&{\rm const.}\,
e^{-\pi({1\over\rho}+{1\over\lambda}-{4\pi\over\rho\lambda n})
|\tilde\gamma_{j,m}-\tilde\gamma_{j,m'}+\beta_{r_{j,m}r_{j,m'}}|},
\\
|g_J(\gamma)|
&\leq&{\rm const.}\,
\prod_{1\leq j\leq n-1}
\Bigl\{
\prod_{1\leq m\leq\nu_j}
e^{-{\pi\over\rho}(\tilde\gamma_{j,m}-\tilde\gamma_{j-1,\mp(J,j,m)})}\no
\\
&\times&\!\!\!\!\!\!\!\!
\prod_{1\leq m\leq\nu_j \atop {1 \le m' \le \nu_{j-1} \atop m'\not=\mp(J,j,m)}}
\!\!\!\!\!\!
e^{{\pi\over\rho}|\tilde\gamma_{j,m}-\tilde\gamma_{j-1,m'}
+\beta_{r_{j,m}r_{j-1,m'}}|}
\!\!\!\!\prod_{1\leq m<m'\leq\nu_j}\!\!\!\!
e^{{\pi\over\rho}
|\tilde\gamma_{j,m}-\tilde\gamma_{j,m'}+\beta_{r_{j,m}r_{j,m'}}|}
\Bigr\},\no
\ena
\bea
|g_{J'}(\gamma_\sigma)| &\leq&{\rm const.}\!\!\prod_{1\leq j\leq n-1}
\Biggl\{\prod_{1\leq m\leq\nu_j}e^{-{\pi\over\lambda}
(\tilde\gamma_{j,\sigma_j(m)}-\tilde\gamma_{j-1,\sigma_{j-1}(\mp(J',j,m))}
+\beta_{r_{j,\sigma_j(m)}r_{j-1,\sigma_{j-1}(\mp(J',j,m))}})}
\no
\ena
\bea
&&\times\prod_{1\leq m\leq\nu_j\atop m'\not=\mp(J',j,m)}\!\!\!\!
e^{{\pi\over\lambda}|\tilde\gamma_{j,\sigma_j(m)}
-\tilde\gamma_{j-1,\sigma_{j-1}(m')}
+\beta_{r_{j,\sigma_j(m)}r_{j-1,\sigma_{j-1}(m')}}|} \\
&&\times\prod_{1\leq m<m'\leq \nu_j}\!\!\!\!
e^{{\pi\over\lambda}|\tilde\gamma_{j,m}-\tilde\gamma_{j,m'}
+\beta_{r_{j,m}r_{j,m'}}|}
\Biggl\}. \no
\ena
Therefore, we have
\bea
|F_{J,J',\sigma}(\gamma)|&\leq&{\rm const.}\,
e^{{2\pi\over\rho\lambda}\sum_{1\leq j\leq n-1\atop1\leq m\leq\nu_j}
(\mu_{j+1}-\mu_j)\tilde\gamma_{j,m}
+{2\pi\over\rho\lambda}\sum_{1\leq r\leq N}\mu_{J_r+1}\beta_r}
\nonumber\\&\times&
\prod_{{1\leq j\leq n-1\atop1\leq m\leq\nu_j}\atop1\leq m'\leq\nu_{j-1}}
e^{-{2\pi^2\over\rho\lambda n}
|\tilde\gamma_{j,m}-\tilde\gamma_{j-1,m'}+\beta_{r_{j,m}r_{j-1,m'}}|}
\nonumber\\&\times&
\prod_{1\leq j\leq n-1\atop1\leq m<m'\leq\nu_j}
e^{{4\pi^2\over\rho\lambda n}
|\tilde\gamma_{j,m}-\tilde\gamma_{j,m'}+\beta_{r_{j,m}r_{j,m'}}|}\label{D}
\\&\times&
\prod_{1\leq j\leq n-1\atop1\leq m\leq \nu_j}
e^{-{\pi\over\rho}\xi(\tilde\gamma_{j,m}-\tilde\gamma_{j-1,\mp(J,j,m)})}
\nonumber\\&\times&
\prod_{1\leq j\leq n-1\atop1\leq m\leq \nu_j}
e^{-{\pi\over\lambda}\xi(\tilde\gamma_{j,\sigma(m)}
-\tilde\gamma_{j-1,\sigma(\mp(J',j,m))}
+\beta_{r_{j,\sigma(m)}r_{j-1,\sigma(\mp(J',j,m))}})}\nonumber
\ena
Here
\bea
\xi(x)=x+|x|.\label{A}
\ena

We apply $-|A+B|\leq|A|-|B|$ to the second line of (\ref{D}),
and $|A+B|\leq|A|+|B|$ to the third line. Then, we use
\bea
\sum_{{1\leq j\leq n-1\atop1\leq m\leq\nu_j}\atop1\leq m'\leq\nu_{j-1}}
|\beta_{r_{j,m}r_{j-1,m'}}|&=&
\sum_{1\leq r<s\leq N}(2\min(J_r,J_s)+1-\delta_{J_r,J_s})
\beta_{sr},\\[12pt]
\sum_{1\leq j\leq n-1\atop1\leq m<m'\leq\nu_j}
|\beta_{r_{j,m}r_{j,m'}}|&=&
\sum_{1\leq r<s\leq N}\min(J_r,J_s)\beta_{sr}.
\ena
We ignore the last line of (\ref{D}). After all these steps, it is enough to
show
\bea
&&{2\pi\over\rho\lambda}\sum_{1\leq j\leq n-1\atop1\leq m\nu_j}
(\mu_{j+1}-\mu_j)\tilde\gamma_{j,m}
+{2\pi^2\over\rho\lambda n}\!\!\!\!
\sum_{{1\leq j\leq n-1\atop1\leq m\leq\nu_j}\atop1\leq m'\leq\nu_{j-1}}
|\tilde\gamma_{j,m}-\tilde\gamma_{j-1,m'}|
\nonumber\\&&
+{4\pi^2\over\rho\lambda n}\!\!\!\!
\sum_{1\leq j\leq n-1\atop1\leq m<m'\leq\nu_j}
|\tilde\gamma_{j,m}-\tilde\gamma_{j,m'}|
-{\pi\over\rho}\sum_{1\leq j\leq n-1\atop1\leq m\leq\nu_j}
\xi(\tilde\gamma_{j,m}-\tilde\gamma_{j-1,\mp(J,j,m)})
\nonumber\\&&\quad
<-\kappa\sum_{1\leq j\leq n-1\atop1\leq m\leq\nu_j}
|\gamma_{j,m}-\gamma_{j-1,\mp(J,j,m)}|.\label{C}
\ena
The left hand side is not larger than
\bea
\sum_{1\leq j\leq n-1\atop1\leq m\leq\nu_j}
\Bigl\{
{2\pi\over\rho\lambda}
(\mu_{\jmax(J,j,m)+1}-\mu_j)(\gamma_{j,m}-\gamma_{j-1,\mp(J,j,m)})
\nonumber\\
+{K\over\rho\lambda}|\gamma_{j,m}-\gamma_{j-1,\mp(J,j,m)}|
-{\pi\over\rho}\xi(\gamma_{j,m}-\gamma_{j-1,\mp(J,j,m)})\Bigr\},
\ena
where
\bea
K&=&{2\pi^2\over n}\!\!\!\!
\sum_{{1\leq j\leq n-1\atop1\leq m\leq\nu_j}\atop1\leq m'\leq\nu_{j-1}}1
+{4\pi^2\over n}\!\!\!\!
\sum_{1\leq j\leq n-1\atop1\leq m<m'\leq\nu_j}1,
\ena
and
\bea
\jmax(J,j,m)&=&\max\{j';r_{j,m}\in{\cal N}_{j'}^{J}\}.
\ena

Choose $\varepsilon,\kappa$ so that
\bea
n\varepsilon+{K\over\rho\lambda}-{2\pi\over\rho}&<&-\kappa,\\
\varepsilon-{K\over\rho\lambda}&>&\kappa.
\ena
This is possible if
\bea
{2K\over\rho\lambda}<{2\pi\over\rho}.
\ena
Then, the estimate (\ref{C}) follows from (\ref{A}) and (\ref{B}).
\qed
\vsk>

For $J \in {\cal Z}_{\nus}$, we set
\bea
P_{J}=\exp{ \Bigl(
\frac{2\pi^2 }{\rho \lambda n}\sum_{1 \le r < s \le N}(1-\delta_{J_{r}, J_{s}})\beta_{sr}
-\frac{2\pi}{\rho \lambda}\sum_{r=1}^{N}\mu_{J_{r}+1}\beta_{r} \Bigr)}.
\label{PJ}
\ena

The following is an obvious consequence of Proposition \ref{prop1}.
\begin{coro}
The integral (\ref{Kdef}) is absolutely convergent. The convergence
is uniform in the variables $\beta$ if we multiply $P_J$ to the integrand.
\end{coro}

Define a partial order in ${\cal Z}_{\nu_1,\ldots,\nu_n}$:
\bea
J\leq J'\quad\hbox{if and only if}\quad J_r+\cdots+J_N\leq J'_r+\cdots+J'_N
\quad\hbox{for all }\ r\,.
\ena

\begin{prop}\label{prop2}
If $J\not\leq J'$, then we have
\bea
\Lim_{\beta_1\ll \cdots\,\ll \beta_N} P_{J}
\left( \prod_{j=1}^{n-1}\prod_{m=1}^{\nu_{j}}
\int_{C_{j}} d\gamma_{j, m} \right)
F_{J,J',\sigma}(\gamma)=0.
\label{E}
\ena
\end{prop}

\proof
We follow the estimate in the proof of Proposition \ref{prop1}.
When we go from (\ref{D}) to (\ref{C}), we dropped the last line
in (\ref{D}). This time we use that term. Namely, we can claim that (\ref{E})
holds unless for some $\sigma=(\sigma_1,\ldots,\sigma_{n-1})$
\bea
r_{j,\sigma_j(m)}^{J}\leq
r_{j-1,\sigma_{j-1}(\mp(J',j,m))}^{J}\label{F}
\ena
holds for all $j$ and $m$. This is clear because
\bea
\xi(x+y)\,=\,2(x+y)\quad{\rm if}\quad y>-x\,.
\ena
We show that (\ref{F}) implies $J\leq J'$. This will complete the proof.

First we prove
\bea
r_{j,\sigma_j(m)}^{J}\leq r_{j,m}^{J'}\quad\hbox{for all }\ m \label{G}
\ena
by induction on $j$. The case $j=0$ is obvious. Suppose that (\ref{G}) is
true for $j-1$. Then we have
\bea
r_{j,\sigma_{j}(m)}^{J}\,\leq\,r_{j-1,\sigma_{j-1}(\mp(J',j,m))}^{J}
\,\leq\,r_{j-1,\mp(J',j,m)}^{J'}\,=\,r_{j,m}^{J'}\,.
\ena
Therefore, (\ref{G}) is true for all $j$. It follows from (\ref{G}) that
\bea
r_{j,m}^{J}\leq r_{j,m}^{J'}\quad\hbox{for all \ $j,m$}\,.
\ena
Finally, we prove $J\leq J'$. This is clear because
\bea
J_r+\cdots+J_N\,=\,\#\{(j,m);r_{j,m}^{J}\geq r.\}
\ena
The proof of Proposition \ref{prop2} is over.\qed
\vsk>

\noindent
This proposition shows that in the asymptotic limit the matrix
$\left(I(w_{J}^{(\rho)}, w_{J'}^{(\lambda)})\right)_{J,J'}$
is triangular.

We have also
\begin{prop}\label{prop3}
\bea
\Lim_{\beta_1\ll \cdots\,\ll \beta_N} P_{J}
\left( \prod_{j=1}^{n-1}\prod_{m=1}^{\nu_{j}}
\int_{C_{j}} d\gamma_{j, m} \right)
F_{J,J,\sigma}(\gamma)=0
\ena
unless $\sigma_j={\rm id}$ for all $j$.
\end{prop}

\proof
Suppose that
\bea
r_{j,\sigma_j(m)}^{J}\leq r_{j-1,\sigma_{j-1}(\mp(J,j,m))}^{J}
\ena
for all $j,m$. From the proof of Proposition \ref{prop2}
we have
\bea
r_{j,\sigma_j(m)}^{J}\leq r_{j,m}^{J}\quad\hbox{for all \ $j,m$.}
\ena
This implies that $\,\sigma_j={\rm id}\;$ for all $j$.
\qed
\vsk>

Define
\bea
\nu^{J, +}_{j,r}=\#\{s\in{\cal N}_{j}^{J};r<s\}\,,
\qquad\nu^{J, -}_{j,r}=\#\{s\in{\cal N}_{j}^{J};r>s\}.
\ena
{}From (\ref{app2}), we have
\begin{prop}\label{prop4}
\bea
&&\Lim_{\beta_1\ll \cdots\,\ll \beta_N}
P_{J}
\left( \prod_{j=1}^{n-1}\prod_{m=1}^{\nu_{j}}
\int_{C_{j}} d\gamma_{j, m} \right)
F_{J,J,{\rm id}}(\gamma)\nonumber\\
&&=
{e^{({1\over\rho}+{1\over\lambda}){\pi^2 i\over n}\sum_{j=1}^{n-1}
\{\nu_j(\nu_{j-1}-1)-\nu_j(\nu_j-1)\}}\over 2^{d_{\nus}}}\prod_{r=1}^N
G_{J_r}(\tilde\mu_{1, r}^{J},\ldots,\tilde\mu_{n, r}^{J})
\label{301}
\ena
where
\bea
G_k(\mu_1,\ldots,\mu_n)&=&
\prod_{j=1}^k\int_{C_{j}}{d\gamma_j\over2\pi i}\prod_{j=1}^k
\varphi(\gamma_j-\gamma_{j-1})
e^{{2\pi\over\rho\lambda}\sum_{j=1}^k(\mu_{j+1}-\mu_j)\gamma_j},\\
\tilde\mu_{j, r}^{J}&=&
\mu_j+{\pi\over n}\sum_{\varepsilon=\pm}\varepsilon
(\nu^{J, \varepsilon}_{j,r}-\nu^{J, \varepsilon}_{j-1,r})
-\delta_{j, J_{r}+1}\frac{\rho+\lambda}{2}.
\ena
In the above formula of $G_{k}$,
$\gamma_{0}=0$ and
the contour $C_{j}$ for $\gamma_{j}$ is a deformation
of the real line $(-\infty, \infty)$ such that the poles at
\bea
\gamma_{j-1}-\frac{\pi i}{n}+\rho i {\Bbb Z}_{\ge 0}+\lambda i {\Bbb Z}_{\ge 0}
\ena
are above $C_{j}$ and the poles at
\bea
\gamma_{j-1}+\frac{\pi i}{n}-\rho i {\Bbb Z}_{\ge 0}-\lambda i {\Bbb Z}_{\ge 0}
\ena
are below $C_{j}$.
\end{prop}
This proposition shows that in the asymptotic limit the diagonal element
$I(w_{J}^{(\rho)}, w_{J'}^{(\lambda)})$ reduces to the one point functions
$G_{J_r}$ $(1\leq r\leq N)$.

Now we consider the asymptotics of $D_{\lambdas}/E_{\lambdas}$.
Hereafter we use the notation $\sim$ as follows:
\bea
f(\beta_{1}, \cdots , \beta_{N}) \sim g(\beta_{1}, \cdots , \beta_{N})
\quad \stackrel{\rm def}{\Longleftrightarrow} \quad
\Lim_{\beta_1\ll \cdots\,\ll \beta_N}\left\{\,
\frac{f(\beta_{1},\cdots, \beta_{N})}{g(\beta_{1},\cdots,\beta_{N})}\,\right\}
\,=\,1\,.
\ena
{}From (\ref{app1}), we have
\bea
E_{\lambdas}(\beta_{1}, \cdots , \beta_{N}) &\sim&
\exp{ \Bigl( \frac{2\pi}{\rho \lambda}
\left( \sum_{j=1}^{n}
{N-1 \atopwithdelims() \lambda_{1}, \cdots , \lambda_{j}-1, \cdots , \lambda_{n}}\mu_{j} \right)
\sum_{m=1}^{N}\beta_{m} \Bigr)} \\
& & {}\times
\exp{ \Bigl( \left\{
\left( \frac{1}{\rho}+\frac{1}{\lambda} \right)\frac{2\pi^2 i}{n}
-\frac{4\pi^{2}}{n\rho \lambda}\sum_{1 \le r<s \le N}\beta_{sr} \right\}
\Lambda_{\lambdas}^{(2)} \Bigr)}. \no
\ena
We note that
\bea
& &\prod_{J \in {\cal Z}_{\nus}} \!\!\!\! P_{J} \\
& &{}=
\exp{\Bigl( \frac{4\pi^{2}}{n\rho \lambda}\sum_{1 \le r<s \le N}\beta_{sr}\Lambda_{\lambdas}^{(2)}
-\frac{2\pi}{\rho\lambda}
\left( \sum_{j=1}^{n}
{N-1 \atopwithdelims() \lambda_{1}, \cdots , \lambda_{j}-1, \cdots , \lambda_{n}}\mu_{j} \right)
\sum_{m=1}^{N}\beta_{m} \Bigr)}. \no
\ena
Hence we find
\bea
& & \frac{D_{\lambdas}(\beta_{1}, \cdots , \beta_{N})}{E_{\lambdas}(\beta_{1}, \cdots , \beta_{N})} \no \\
& & {}\sim
\exp{\Bigl(-\left(\frac{1}{\rho}+\frac{1}{\lambda}\right)\frac{2\pi^{2}i}{n}
\Lambda_{\lambdas}^{(2)} \Bigr)}
\det{\left(P_{J} I(w_{J}^{(\rho)}, w_{J'}^{(\lambda)})
\right)_{J, J' \in {\cal Z}_{\nus}}}.
\ena
{}From Propositions \ref{prop2}, \ref{prop3}, and \ref{prop4}, we see that
\bea
\det{\left(P_{J} I(w_{J}^{(\rho)},w_{J'}^{(\lambda)})
\right)_{J,J'\in {\cal Z}_{\nus}}}\sim
\prod_{J\in{\cal Z}_{\nus}}(\hbox{the right hand side of \,(\ref{301})}).
\ena
Therefore, we get
\begin{prop}\label{prop4.2}
\bea
\kern-4em &&
c_{\lambdas}(\mu_{1}, \cdots , \mu_{n} ; \rho , \lambda )\,=\,
2^{-d_{\nus}\Lambda_{\lambdas}^{(0)}}
\label{203}
\\[6pt]
\kern-4em &&
{}\times\,\exp{\Bigl((N^{2}-N-2)\left(\frac{1}{\rho}+\frac{1}{\lambda}
\right) \frac{\pi^{2}i}{n}\;\Lambda_{\lambdas}^{(2)} \Bigr)}\,
\prod_{J \in{\cal Z}_{\nus}}\prod_{r=1}^{N}G_{J_{r}}
(\widetilde{\mu}_{1, r}^{J}, \cdots , \widetilde{\mu}_{n, r}^{J})\,.\no
\ena
\end{prop}

\proof
Note that
\vv-.6>
\bea
\# {\cal Z}_{\nus}=\Lambda_{\lambdas}^{(0)}.
\ena
We get the term
\vv-.5>
\bea
\exp{\Bigl( (N^{2}-N-2)\left( \frac{1}{\rho}+\frac{1}{\lambda}\right)
\frac{\pi^{2}i}{n}\;\Lambda_{\lambdas}^{(2)}\Bigr)}
\ena
by using the following formulae:
\bea
& & \sum_{j=1}^{n-1}\{ \nu_{j}(\nu_{j-1}-1)-\nu_{j}(\nu_{j}-1) \}=
\frac{1}{2}\left( N^{2}-\sum_{j=1}^{n}\lambda_{j}^{2} \right),
\no \\
& & \frac{1}{2}\left(N^{2}-\sum_{j=1}^{n}\lambda_{j}^{2} \right)\Lambda_{\lambdas}^{(0)}=
N(N-1)\Lambda_{\lambdas}^{(2)}.
\ena
$\square$ % \qed

%\section{The $N=1$ case}
\section{Proof of Theorem \protect\ref{thm6.1}}

%In this section
First, we find an explicit formula for
$G_{k}(\mu_{1}, \cdots , \mu_{n})$. We set %!
\bea
H_{k}(x_{1}, \cdots , x_{k})=\left( \prod_{j=1}^{k}\int_{C_{j}}d\gamma_{j} \right)
\prod_{j=1}^{k}\varphi(\gamma_{j}-\gamma_{j-1})
e^{\frac{2\pi}{\rho\lambda}\sum_{j=1}^{k}x_{j}(\gamma_{j}-\gamma_{j-1})},
\label{500}
\ena
where $\gamma_{0}=0$. Then we have
\bea
G_{k}(\mu_{1}, \cdots , \mu_{n})=H_{k}(\mu_{k+1}-\mu_{1}, \cdots,
\mu_{k+1}-\mu_{k}).
\label{501}
\\[-.9\bls]\no
\ena
The integral (\ref{500}) is absolutely convergent if
\bea
|\Re x_{j}|\,<\,\frac{\rho+\lambda}{2}+\frac{\pi}{n},\qquad(j=1,\cdots,n).
\ena
By changing the integration variables $\gamma_{j}$ to
\vv.3>
\bea
u_{j}=\gamma_{j}-\gamma_{j-1}, \qquad (j=1, \cdots , n)\,,
\\[-.7\bls]\no
\ena
we can see that
\vv->
\bea
H_{k}(x_{1}, \cdots , x_{k})=\prod_{j=1}^{k}H(x_{j}),
\label{relH}
\ena
where
\vv->
\bea
H(x)=\int_{C}du\, \varphi(u)\, e^{\frac{2\pi}{\rho\lambda}xu}.
\label{defH}
\ena
In the above formula, the contour $C$ is a deformation of the real line $(-\infty, \infty)$
such that the poles at
\bea
-\frac{\pi i}{n}+\rho i{\Bbb Z}_{\ge 0}+\lambda i {\Bbb Z}_{\ge 0}
\ena
are above $C$ and the poles at
\bea
\frac{\pi i}{n}-\rho i {\Bbb Z}_{\ge 0}-\lambda i {\Bbb Z}_{\ge 0}
\ena
are below $C$.

The explicit formula for the function $H$ is obtained in \cite{MT1}.

\begin{prop}\label{5.1}
\bea
H(x)=\frac{\sqrt{\rho\lambda}}{S_{2}(-\frac{2\pi}{n})}
     \frac{S_{2}(x+\frac{\rho+\lambda}{2}-\frac{\pi}{n})}{S_{2}(x+\frac{\rho+\lambda}{2}+\frac{\pi}{n})}
\label{form1}
\ena
\end{prop}

\nt
{}From (\ref{501}), (\ref{relH}) and (\ref{form1}), we get

\begin{prop}\label{prop5.4}
\bea
G_{k}(\mu_{1}, \cdots , \mu_{n})=
\left( \frac{\sqrt{\rho\lambda}}{S_{2}(-\frac{2\pi}{n})} \right)^{k}
\prod_{j=1}^{k}\frac{S_{2}(\mu_{k+1}-\mu_{j}+\frac{\rho+\lambda}{2}-
\frac{\pi}{n})}
{S_{2}(\mu_{k+1}-\mu_{j}+\frac{\rho+\lambda}{2}+\frac{\pi}{n})}\;.
\label{512}
\ena
\end{prop}

%\section{Proof of Theorem \protect\ref{thm6.1}}

\nt
Now it remains to calculate %!
\bea
\no\\[-24pt]
\prod_{J \in {\cal Z}_{\nus}} \prod_{r=1}^{N}
G_{J_{r}}(\widetilde{\mu}_{1, r}^{J}, \cdots , \widetilde{\mu}_{n, r}^{J}).
\ena
We set %!
\bea
M_{j, k}^{J, +}= \# \{ r ; J_{r}=j, k<r \}, \quad
M_{j, k}^{J, -}= \# \{ r ; J_{r}=j, k>r \}.
\ena
Note that
\bea
\sum_{r=1}^{n}J_{r}=\sum_{j=1}^{n}(j-1)\lambda_{j}, \quad {\rm for \,all} \quad
J \in {\cal Z}_{\nus}.
\ena
{}From (\ref{301}) and (\ref{512}), we have
\bea
& &
\prod_{J \in {\cal Z}_{\nus}} \prod_{r=1}^{N}
G_{J_{r}}(\widetilde{\mu}_{1, r}^{J}, \cdots , \widetilde{\mu}_{n, r}^{J}) \\
& & {}=
\left( \frac{\sqrt{\rho\lambda}}{S_{2}(-\frac{2\pi}{n})} \right)
^{\Lambda_{\lambdas}^{(0)}\sum_{j=1}^{n}(j-1)\lambda_{j}}\!\!\!\!\!\!
\prod_{1 \le r'<r \le n}\prod_{J \in {\cal Z}_{\nus}}\prod_{k \atop J_{k}+1=r}
\frac{S_{2}(\mu_{r}-\mu_{r'}+\frac{\pi}{n}(D_{r', r, k}^{J}-1))}
{S_{2}(\mu_{r}-\mu_{r'}+\frac{\pi}{n}(D_{r', r, k}^{J}+1))}, \no
\ena
where $D_{r', r, k}^{J}$ is given by
\bea
D_{r', r, k}^{J}=\sum_{\eps=\pm}
\eps (M_{r'-1, k}^{J, \eps}-M_{r-1, k}^{J, \eps})
=\lambda_{r'}-\lambda_{r}+1-2(M_{r'-1, k}^{J, -}-M_{r-1, k}^{J, -})
\ena
for $k$ satisfying $J_{k}+1=r$.

Now we rewrite
\bea
\prod_{J \in {\cal Z}_{\nus}}\prod_{k \atop J_{k}+1=r}
\frac{S_{2}(\mu_{r}-\mu_{r'}+\frac{\pi}{n}(D_{r', r, k}^{J}-1))}
{S_{2}(\mu_{r}-\mu_{r'}+\frac{\pi}{n}(D_{r', r, k}^{J}+1))}\;.
\label{601}
\ena
Let us consider the following set:
\bea
{\cal F}^{(r', r)}=
\bigsqcup_{(J, k) \atop {J \in {\cal Z}_{\nus} \atop J_{k}+1=r}}\{M_{r'-1, k}^{J, -}-M_{r-1, k}^{J, -}\},
\quad (1 \le r' < r \le n),
\ena
where $\bigsqcup$ means a disjoint union.
For $a \in {\bf Z}$, we set
\bea
{\rm mult}^{(r', r)}(a)=\# \{ t \in {\cal F}^{(r', r)}; t=a \}.
\ena
Then we have
\bea
(\ref{601})=\prod_{a \in {\bf Z}}S_{2}\left(
\mu_{r}-\mu_{r'}+\frac{\pi}{n}(\lambda_{r'}-\lambda_{r}-2a)\right)
^{{\rm mult}^{(r', r)}(a)-{\rm mult}^{(r', r)}(a+1)}.
\ena
We can show
\bea
& & {\rm mult}^{(r', r)}(a)-{\rm mult}^{(r', r)}(a+1) \no \\[6pt]
& & {}=\left\{
\begin{array}{ll}
\displaystyle{
-\Lambda_{\lambdas}^{(0)}
{\lambda_{r}+\lambda_{r'} \atopwithdelims() \lambda_{r}+a}\bigg/
{\lambda_{r}+\lambda_{r'} \atopwithdelims() \lambda_{r}} }\,,&
-\lambda_{r} \le a \le \min\{-1, \lambda_{r}-\lambda_{r'}\}, \\[12pt]
\displaystyle{
\Lambda_{\lambdas}^{(0)}
{\lambda_{r}+\lambda_{r'} \atopwithdelims() \lambda_{r'}-a}\bigg/
{\lambda_{r}+\lambda_{r'} \atopwithdelims() \lambda_{r}} }\,,&
\max\{0, \lambda_{r'}-\lambda_{r}+1\} \le a \le \lambda_{r'}, \\[12pt]
0\,, &{\rm otherwise}.
\end{array}
\right. \no
\ena
This completes the proof. \qed

\addtocounter{section}{1}
\setcounter{equation}{0}

\section*{Appendix}

Here we summarize the property of the double sine function
$S_2(x)=S_{2}(x|\om_{1},\om_{2})$ following \cite{JM}.

We assume that $\Re\om_{1} >0, \Re\om_{2}>0$.
$S_{2}(x|\om_{1},\om_{2})$ is a meromorpic function of $x$ and
symmetric with respect to $\om_{1},\om_{2}$. Its zeros and poles are
given by
$$
\hbox{zeros at } \ x=\om_{1}{\Bbb Z}_{\le 0}+\om_{2}{\Bbb Z}_{\le 0}\,,\qquad
\hbox{poles at } \ x=\om_{1}{\Bbb Z}_{\ge 1}+\om_{2}{\Bbb Z}_{\ge 1}\,.
$$
Its asymptotic behavior is as follows:
\bea
\label{app1}
\log S_{2}(x)\,=\,\pm \pi i \! \left( \frac{x^2}{2 \om_{1} \om_{2}}-
\frac{ \om_{1}+
\om_{2}}{2 \om_{1} \om_{2}}\,x-\frac{1}{12}\left(\frac{\om_{1}}{\om_{2}}+
\frac{\om_{2}}{\om_{1}}+3 \right) \right)\! \! + \! o(1)\,, &&
\\[8pt]
(x \to \infty, \,\pm\Im x > 0)\,. && \no
\\[-2\bls] \no
\ena
This implies that
\bea
\log S_{2}(a+x)S_{2}(a-x)\,=\,
\pm \pi i \frac{2a-\om_{1}-\om_{2}}{\om_{1} \om_{2}}\,x+o(1),
\quad (x \to \infty\,,\;\ \pm \Im x >0)\,. \label{app2}
\ena
The following formulae hold:
\bea
\frac{S_{2}(x+\om_{1})}{S_{2}(x)} &=& \frac{1}{2 \sin \frac{\pi x}{\om_{2}}},
\label{app3} \\
S_{2}(x) &=& \frac{2 \pi}{\sqrt{\om_{1}\om_{2}}}x+O(x^2) \quad (x \to 0)\,.
\label{app5}
\ena

\ifUS\else\newpage\fi

\end{document}